\documentclass[11pt]{amsart}
\hoffset= - 0.75 in

\usepackage{amsfonts, amsmath, amssymb, amsthm}
\usepackage{latexsym}
\newtheorem{thm}{Theorem}
\newtheorem{claim}{Claim}
\newtheorem{conj}{Conjecture}
\newtheorem{obs}{Observation}
\newtheorem{lem}{Lemma}
\newtheorem{cor}{Corollary}
\newtheorem{prop}{Proposition}
\newtheorem{fact}{Fact}
\newcommand{\khalf}{\frac{k}{2}}
\newcommand{\ZZ}{\mathbb{Z}}
\newcommand{\NN}{\mathbb{N}}
\newcommand{\Zn}{\ZZ/n\ZZ}
\newcommand{\GU}{G_U}
\newcommand{\Gob}{G_{\{0,b\}}}
\newcommand{\bin}[2]{
	\left(
		\begin{array}{@{}c@{}}
			#1  \\  #2
		\end{array}
	\right)		}

\begin{document}

\title{Avoidable Sets in Groups}
\author{Mike Develin}
\address{Department of Mathematics, University of California-Berkeley,
Berkeley, CA, 94720-3840}
\curraddr{2706B Martin Luther King, Berkeley, CA 94703-2110}
\date{\today}
\email{develin@math.berkeley.edu}

\begin{abstract}
In a set equipped with a binary operation, $(S,\cdot)$, a
subset $U$ is defined to be avoidable if there exists a partition
$\{A,B\}$ of $S$ such that no element of $U$ is the product of two
distinct elements of $A$ or of two distinct elements of $B$.
For more than two decades, avoidable sets
in the natural numbers (under addition) have been studied by
renowned mathematicians such as Erd\"{o}s, and a few families of 
sets have been shown to be avoidable in that setting. In this paper we 
investigate the generalized notion of an avoidable
set and determine the avoidable sets in several families of groups;
previous work in this field considered only the case $(S,\cdot)=(\NN,+)$.
\end{abstract}

\maketitle

\section{Introduction}
In a set equipped with a binary operation, $(S,\cdot)$, a
subset $U$ is defined to be \textbf{avoidable} if there exists a partition
$\{A,B\}$ of $S$ such that no element of $U$ is the product of two
distinct elements of $A$ or of two distinct elements of $B$; the partition
$\{A,B\}$ is said to \textbf{avoid} $U$.

Avoidable sets were first studied by Alladi, Erd\"{o}s, and Hoggatt
(\cite{Erdos}) in 1978. However, throughout the intervening decades,
avoidable sets have been studied only in the case
$(S,\cdot)=(\mathbb{N},+)$. All of the 
relevant concepts in this particular case generalize in the context of the
definition
presented above: the operation need not be 
commutative or even associative. In the case of $(\mathbb{N},+)$, only a few 
families of sets (\cite{Chow}, \cite{Shan}, \cite{Develin}) have been shown 
to be avoidable, and almost no general facts are known 
about the nature of avoidable sets. In this paper, we discuss 
avoidable sets in a different setting: namely, when $S$ is a group.

\section{Abelian Groups}
Note: In this section, all groups will be written additively, and the 
identity element will be denoted by 0.

First, we present two basic definitions that will be used throughout
the paper. We say that $U$ is 
\textbf{saturated} if $U$ is
maximal among the collection of avoidable subsets of $S$. Also, for any
subset $U$, 
we define the 
\textbf{associated graph}, $\mathbf{G_U}$, to be the graph whose
vertices are the elements of $S$, where two elements are connected by an edge
if and only if their sum is an element of $U$; then 
$U$ is avoidable if and only if $G_U$ is 2-colorable. 

Abelian groups present a particularly nice case to study: if $|U|=n$, then
each vertex $x$ in $G_U$ will have degree $n$ or $n-1$, depending
on whether $2x$ is in $U$ or not. In this section we discuss the avoidable
sets in abelian groups, commencing with a complete categorization of the
avoidable sets in $\ZZ/n\ZZ$, and proceeding with a discussion on the
nature of avoidable sets in arbitrary finitely generated abelian groups.
Throughout this section, we will call an element $b$ \textbf{even} if there exists
$x$ such that $2x=b$, and \textbf{odd} otherwise. 

\subsection{Saturated Sets in $\ZZ/n\ZZ$}\label{Zn}

In this section, all equalities, unless stated otherwise, are in $\Zn$.
 
We first note that since all avoidable sets in $\Zn$ are finite, to categorize all
avoidable sets we need merely categorize the saturated sets.  We present first 
two lemmas about unavoidable sets.

\begin{lem}
\label{4lem} No four-element set is avoidable in $\Zn$
for n=2k, unless it is composed solely of odd numbers.
\end{lem}
\noindent

\begin{proof}
Consider a four-element set $\{a,b,c,d\}$, with $a$ even. Two of $b$, $c$,
and $d$ have
the same parity, without loss of generality $b$ and $c$. But then
if we pick representatives for the residue classes $a$, $b$, and $c$
in $\ZZ$, $\frac{a+b-c}{2}$, $\frac{a-b+c}{2}$, and $\frac{-a+b+c}{2}$
are all integers, and the corresponding residue classes (distinct as
$a$, $b$, and $c$ are) form a 3-cycle in $G_{\{a,b,c,d\}}$.
\end{proof}

\begin{lem}\label{5lem}
In an arbitrary abelian group, $\{0,b,c\}$ is not avoidable unless either 
$2b=c$, $2c=b$, or $2b=2c$.
\end{lem}
\noindent
\begin{proof}
Unless one of these conditions holds, $0, b, c-b, b-c, c$ 
is a 5-cycle in $G_{\{0,b,c\}}$. Note that if one of these conditions 
holds, then two neighboring elements of this supposed 5-cycle are 
identical, and so the 5-cycle does not exist.
\end{proof}

We now state the main theorem proven in this section.

\begin{thm}\label{cyclic}
The saturated sets in $\ZZ/n\ZZ$ are:

(a) if $n=2k+1$, $\{a,b\}$ with $a\ne b$

(b) if $n=4k+2$, $\{a,a+x,a+2x\}$ with $x$ odd, $x\ne n/2$, and $a$ even;
$\{a,a+n/2\}$; and $\{x\,|\,x=2l+1, 1\le l\le n/2\}$ if $k>0$

(c) if $n=4k$, $\{a,a+x,a+2x\}$ with $x$ odd, $a$ even; $\{a,a+4x\}$ with $a$
even, $x\ne n/4$; and $\{x\,|\,x=2l+1, 1\le l\le n/2\}$
\end{thm}
\noindent
\begin{proof}
(a) In this case, as 2 and $n$ are relatively prime, 2 has a 
multiplicative inverse in  $\Zn$. Denote this element by $y$. 

First, we note that no three-element set is avoidable: 
$G_{\{a,b,c\}}$ contains the 3-cycle $(y(a+b-c), y(a-b+c),
y(-a+b+c))$; these elements are distinct, as $a$, $b$, and $c$ are.
On the other hand, any two-element set $\{a,b\}$ is avoidable. To prove
this, it suffices to show that any set of the form $\{0,t\}$ is
avoidable, for then we can construct a partition $\{A,B\}$ avoiding
$\{0,b-a\}$, and add $ya$ to each member of the partition to obtain a
partition avoiding $\{a,b\}$. 

Consider $G_{\{0,t\}}$. Each vertex has degree at most 2,
so its connected components are simply paths and cycles.
Then $G_{\{0,t\}}$ is 2-colorable if and only if all of the
cycles are even. Let us consider an arbitrary cycle, $(x, t-x, x-t,
2t-x, \ldots)$. Because this is a cycle, all of its vertices have
degree 2 and thus no two adjacent terms in this sequence are identical.
If the cycle had odd length, we would have $x=lt-x$ for some $l\in\ZZ$;
however, this implies either $x-mt=mt-x$, if $l=2m$, or
$x-mt=(m+1)t-x$, if $l=2m+1$. In either case, two adjacent terms are
identical, contrary to the assumption that the terms are contained in a cycle.

(b) Note that in cases (b) and (c), the set of all odd numbers is 
avoidable: simply consider the partition $A=\{2l+1\}$, $B=\{2l\}$, 
which obviously avoids the set of all odd numbers. Clearly 
this is also the unique partition avoiding it, and for $n>2$, this set is 
therefore saturated.

Because of Lemma~\ref{4lem}, all that remains is to determine when a 
three-element 
set $\{a,b,c\}$ is avoidable, where $a$ is even. As in (a), we need 
only consider the case $a=0$, because the existence of a partition 
avoiding $\{a,b,c\}$ is equivalent to the existence of a partition 
avoiding $\{0,b-a,c-a\}$.
Also, note that if $a+b+c$ is even, then as before 
$\{a,b,c\}$ cannot be avoidable. So we need only consider the case of 
$\{0,b,c\}$ where $b$ is even and $c$ is odd. Our strategy will be to 
fix $b$ and determine for what values of $c$ the set $\{0,b,c\}$ is 
avoidable.

Denote by $\frac{b}{2}$ the unique even element satisfying 
$2\frac{b}{2}=b$. Then $\Gob$ has four vertices of degree 
1: 0, $\frac{b}{2}$, $2k+1$, and $2k+1+\frac{b}{2}$. Consequently, it is 
composed of two paths and some number of cycles. Within a given 
connected component, all the elements have the same parity, since both 0
and  $b$ are even. Therefore, the path starting at 0 must end at
$\frac{b}{2}$. This path 
consists of $\{0, b, -b, 2b, \ldots, lb=\frac{b}{2}\}$ for some 
$l\in\ZZ$. The last condition 
is equivalent to $(2l-1)b=0$, meaning that the solution first encountered 
is when $l=\frac{1}{2}(1-\frac{n}{(n,b)})$ where $(n,b)$ is the greatest 
common divisor of $n$ and $b$. Note that as 2 appears
with 
multiplicity 1 in both $n$ and $b$, $\frac{n}{(n,b)}$ is odd; therefore,
$l$ is an integer. The 
length of the path is $\frac{n}{(n,b)}$. It is not hard to see that the 
other path is a translate of this path by $2k+1$.

Now let us consider the cycles: an arbitrary cycle is of the form $\{x,
b-x, x-b, 2b-x, \ldots, x-lb\}$. 
Solving for minimal $l$ yields $l=\frac{n}{(n,b)}$, whence the length of 
the cycle is $\frac{2n}{(n,b)}$. So $\{0,b\}$ is certainly avoidable.

We now apply Lemma \ref{5lem}. Since $2b=c$ is 
impossible (as $c$ is odd), and $2b=2c$ is impossible (for then 
$c=b+(2k+1)$, 
whence $G_{\{0,b,c\}}$ contains the proper 5-cycle 
$\frac{b}{2}, -\frac{b}{2}, (2k+1)+\frac{3b}{2}, (2k+1)-\frac{b}{2}, 
(2k+1)+\frac{b}{2}$), for $\{0,b,c\}$ to be avoidable we need $2c=b$. 
As $c$ is odd, this means $c=(2k+1)+\frac{b}{2}$.

But now $\{0,b,c\}$ is avoidable for this value of $c$. To prove this, it 
suffices to construct a \textbf{c-consistent} 2-coloring of $\Gob$, i.e.
a 
2-coloring of $\Gob$ such that $x$ and $c-x$ always have opposite colors. 
Using the term \textbf{c-complement of $x$} to denote $c-x$, we find that:

\begin{fact}\label{parity}
The $c$-complement of $b-x$ is the 0-complement of $c-x$, and the $c$-complement 
of $-x$ is the $b$-complement of $c-x$.
\end{fact}

Fact~\ref{parity}, along with the fact that the $c$-complement of $0$
is $c$, implies that the $c$-complements of a path in
$\Gob$
are completely contained in the other path. We color 
one of the two paths arbitrarily; this induces a consistent coloring on
the other path.
Similarly, the $c$-complements of a given cycle in $\Gob$ lie completely in 
some cycle. This second cycle is distinct from the first; if not,
we would have either $x-lb=c-x$ or $lb-x=c-x$, each
impossible due to parity considerations.

As all cycles have the same length, these complements comprise the
entirety of the second cycle. We can then pair off the cycles
into pairs which are
$c$-complements of each other and color one cycle of each pair arbitrarily;
this induces a consistent coloring on the other cycle.  Therefore, sets of
the form $\{0,b,\frac{b}{2}+(2k+1)\}$ (or equivalently, sets of the form
$\{0,c,2c\}$ for $c$ odd and $c\ne \frac{n}{2}$) are avoidable, and
saturated by Lemma~\ref{4lem}. The only two-element set containing 0
not contained in any of these sets is $\{0,2k+1\}$, which is therefore
saturated. 

(c) The proof given for part (b) works when $b\equiv 2$ (mod 4), with the 
minor modification that now either value of $x$ for which $2x=b$ yields 
an avoidable set $\{0,b,x\}$. Meanwhile, if $b\equiv 0$ (mod 4), 
both such values of $x$ are even, and consequently there are no 3-element 
avoidable sets containing $\{0,b\}$, so $\{0,b\}$ is saturated.

This completes the proof of Theorem~\ref{cyclic}.
\end{proof}

\subsection{Saturated Sets in $\ZZ$}

The case of $\ZZ$ is similar to that of $\Zn$; the proof of 
Lemma~\ref{4lem} carries over, so that no four-element sets are 
avoidable, and the set of all odd numbers is saturated as before. We 
then have the following result.

\begin{thm} The saturated sets in $\ZZ$ are $\{x\,|\,x=2l+1, l\in \ZZ\}$;
$\{a,a+x,a+2x\}$ with $a$ even, $x$ odd; and $\{a,a+4x\}$ with $a$
even. 
\end{thm} 

\begin{proof} We assume that the set contains an even
element; otherwise, it is contained in the first set above. As before,
we can reduce to determining when $\{0,b,c\}$ is avoidable where $b$ is
even and $c$ odd. Due to Lemma~\ref{5lem} we need only consider the
case $c=\frac{b}{2}$. In this case, we look at $\Gob$. There are no
cycles: the paths containing a vertex of degree 1 are $0,b,-b,2b,
\ldots$ and $\frac{b}{2}, -\frac{b}{2}, \frac{3b}{2}, -\frac{3b}{2},
\ldots$, which contain every element of $\ZZ$ which is congruent to 0
or $\frac{b}{2}$ (mod $b$). The other paths are of the form $\ldots,
-b-x, x+b, -x, x, b-x, x-b, 2b-x, 2x-b, \ldots$, and contain every
element congruent to $x$ or $-x$ (mod $b$). Fact~\ref{parity} is still
true, and because these paths contain every element in a given residue
class and its inverse, the $c$-complements of one path still comprise the
entirety of another path. Therefore, we can $c$-consistently color the
paths as in the proof of Theorem~\ref{cyclic}, so that $\{0,b,c\}$ is
avoidable. 

As before, all two-element sets are avoidable; the ones not contained in
any three-element avoidable set are of the form $\{a,a+4x\}$ with $a$
even, and these are consequently saturated. Thus, the proof of the theorem
is complete. 
\end{proof}

\subsection{Finitely Generated Abelian Groups}
In the case of finitely generated abelian groups, we have the following theorem.
\begin{thm}\label{abgroup}
Let $G$ be any finitely generated abelian group, and write $G$ as:

\[
G=\left(\bigoplus_{i=1}^{m_1} \frac{\ZZ}{(2n_i+1)\ZZ}\right)\oplus\left( 
\bigoplus_{i=1}^{m_2} 
\frac{\ZZ}{2\ZZ}\right)\oplus\left( \bigoplus_{i=1}^{m_3} 
\frac{\ZZ}{(4q_i)\ZZ}\right)\oplus\left( \bigoplus_{i=1}^{r} \ZZ\right)\]

Then the maximal size, $m$, of an avoidable set in $G$ containing an even 
number is:

a) $1+2^{m_2-1}$, if $m_1=m_3=r=0$ 

b) $2$, if $m_2=m_3=r=0$

c) $2+2^{m_2-1}$, if $m_3=r=0$ and $m_2\ne 0$

d) $2+2^{m_2}$ otherwise.
\end{thm}
\begin{proof}
Because of previous reasoning, it suffices to find the maximum size of
an avoidable set containing 0; in what follows, all avoidable sets
mentioned are assumed to contain 0. We break the proof into cases
corresponding to the cases of the theorem.

(a) $m_1=m_3=r=0$. In this case, 
$G=(\frac{\ZZ}{2\ZZ})^{m_2}$, and we want to show 
that $m=1+2^{m_2-1}$. Certainly $m\le 1+2^{m_2-1}$, for in any partition one 
of the partition's sets must contain at least $\frac{|G|}{2}$ elements, 
producing at least $\frac{|G|}{2}-1=2^{m_2-1}-1$ distinct sums, so at 
most $|G|-(2^{m_2-1}-1)=1+2^{m_2-1}$ sums are avoidable. However, there 
is a set of cardinality $1+2^{m_2-1}$ which is avoidable, namely $\{0\} 
\cup \{(1 ,x_2, \ldots, x_{m_2})\}$; the partition is the latter 
set and its complement.

(b) $m_2=m_3=r=0$. This case is also simple. Here, 
$G={\displaystyle\bigoplus_{i=1}^{m_1}}\frac{\ZZ}{(2n_i+1)\ZZ}$ and so 
$2$ has a multiplicative inverse $y_i$ in each summand. Using the logic
employed in Theorem~\ref{cyclic}, we conclude that all two-element sets
and no three-element sets are avoidable, so $m=2$.

(c) $m_3=r=0$ and $m_2\ne 0$. Define $S_b=\{0,b\}\cup\{x\,|\,2x=b\}$. By 
Lemma~\ref{5lem}, 
every avoidable set must be a subset of $S_b$
for some $b$. We seek to maximize the size of such a set. If $b$ is 
odd, $|S_b|=2$, which is not of much help, so we assume $b$ is even. Because of
the following lemma, we can assume $b$ is in the set.

\begin{lem}\label{bisgood}
If $U$ is an avoidable subset of $S_b$ with $|U|\ge 2$ and $b\notin U$, then $U 
\cup \{b\}$ is also avoidable. \end{lem}

\begin{proof}
Let $\{A,B\}$ be a partition avoiding $U$. Suppose there are two elements $u$ 
and $v$ 
in one of the sets, without loss of generality $A$, whose sum 
is $b$. Then since $U$ contains some $x$ for which $2x=b$, we know that 
$x-u$ and $x-v$ are both in $B$. But then $(x-u)+(x-v)=2x-(u+v)=0$, 
contradicting the fact that 0 is contained in $U$. Therefore, $\{A,B\}$ 
also avoids $b$, so that $U\cup \{b\}$ is also avoidable.
\end{proof}

We show first that if $b=0$, the maximal size of an avoidable subset $U$ 
of $S_b$ is within the bound in question. In this case, we 
have $|S_0|=2^{m_2}$.

Take some $x\in U$ which is to be in our avoidable set along with 0, and
consider the associated graph $G_{\{0,x\}}$. The induced subgraph on $S_0$
will consist of disjoint copies
of $K_2$, with 0 connected to $x$. However, if $y,z\in S_0$ are connected
by an edge in $G_{\{0,x\}}$, $y$ and $z$ cannot both be in $U$, as then
$0,y,z$ would be a 3-cycle in $G_U$. Therefore, $U$ has at most
$1+2^{m_2-1}$ elements.

Now, we consider the case where $b\ne 0$. In this case, the connected
components of $\Gob$ are paths and cycles; the vertices of degree 1 are
$\{x | 2x=b\}\cup \{x\,|\,2x=0\}$, a set with $2^{m_2+1}$ 
elements.
Therefore, there are $2^{m_2}$ paths; the path
starting at 0 is $0, b, -b, 2b, \ldots, kb=x$, where $2x=0$ or
$2x=b$.
In the first case, we have $2kb=0$; in the second, $2kb=b$.  Looking at
each
coordinate separately, we see that if $k$ is a minimal solution to the
first equation, $k$ must be odd (otherwise $\frac{k}{2}$ also satisfies
the first equation), but then $-\frac{k+1}{2}$ satisfies the second
equation. So a solution to the second equation happens first, and thus the
path starting at 0 ends at a value of $x$ for which $2x=b$, and also has
an even number of edges as this solution is 
negative. All the other paths are
of the form $y, b-y, y-b, 2b-y, \ldots, kb-y=x-y$ where $2y=0$ (and
consequently $2(x-y)=b$).

Let us now pick another element $c\in U$, $2c=b$, and consider the
resulting associated graph $G_{\{0,b,c\}}$. 
Because of Fact~\ref{parity}, the components which contain the vertices 
of the paths in $\Gob$ are now pairs of entangled paths; Each
path contains 
two values of $x$ such that $2x=0$, and two values of $y$ such that 
$2y=b$. These values must satisfy 
$x_1+y_2=x_2+y_1=c$, where $x_1$ and $y_1$ are the endpoints of a path in 
$\Gob$ as are $x_2$ and $y_2$.  Therefore, $x_1$ and $y_1$ 
must be colored identically (because the path has even length), and 
$x_1$ and $y_2$ must be colored differently (because their sum is $c$), 
so $y_1$ and $y_2$ must be colored differently, and hence must be in
opposite sets of any partition avoiding $U$. Consequently, at most one of
$y_1$ and $y_2$ can be in $U$, 
so at most half of the $2^{m_2}$
values of $y$ with $2y=b$ can be in our avoidable set. These values are 
in addition to 0 and $b$, so our avoidable set can have at most 
$2+2^{m_2-1}$ elements. To complete the proof of case (c), it remains only
to show that this bound can be achieved.
We prove this by induction on $m_2$ as follows.

In the base case, $m_2=1$, we wish to find an avoidable set of size 3.
However, this is not hard: the set
$\{(0,\ldots,0), (1,\ldots,1), (2,\ldots,2)\}$ is easily seen to be
avoidable as in Theorem~\ref{cyclic} (b). 

Now, assume that in some abelian group $G$ we have an avoidable set $S$ of
size $2+2^n$ consisting of $\{0,b\}\cup T$ where $T$ is a set of elements
$x$ such that $2x=b$. We need only to show that there is an avoidable
set of size $2+2^{n+1}$ in $G^{\prime}=G\oplus\frac{\ZZ}{2\ZZ}$. Let 
$\{C,D\}$ be a partition of $G$ avoiding $S$. Then we let:
\[
C^{\prime}=
\{(c,s)\,|\,c\in C, s\in \frac{\ZZ}{2\ZZ}\}
\]
and \[
D^{\prime}=\
\{(d,s)\,|\,d\in D, s\in \frac{\ZZ}{2\ZZ}\}
\]

Then $\{C^{\prime}, D^{\prime}\}$ is clearly a partition of $G^{\prime}$
avoiding the $2+2^{n+1}$-element set $S^{\prime}=\{(0,0),
(b,0)\}\cup \{(t,s)\,|\,t\in T, s\in\frac{\ZZ}{2\ZZ}\}$, as desired.

(d) The same construction as in part (c), with the base case replaced by
$m_2=0$, ensures that there exists some avoidable set of size $2+2^{m_2}$.
So we need only prove that 
$2+2^{m_2}$ is an upper 
bound on the size of an avoidable set in $G$. As in part (c), every 
avoidable set $U$ 
is a subset of $S_b$ for some $b$. Consider the set $\{x\,|\,2x=b\}$.
No 
two elements $y$ and $z$ which are equal in all $\frac{\ZZ}{2\ZZ}$ 
summands of $G$ 
can both be in $U$, as then we would have 
$0+y+z$ even, impossible by previous logic.
There are precisely $2^{m_2}$ possible values an element of this set 
can take in 
${\displaystyle\bigoplus_{i=1}^{m_2}}\frac{\ZZ}{2\ZZ}$, so at most 
$2^{m_2}$ of these elements can be in $U$, and thus $U$ can have a total 
of at most $2+2^{m_2}$ elements in all.
\end{proof}

We have found the maximum size of a saturated set containing 0 in an
arbitrary finitely generated abelian group. Indeed, using the methods of
the proof of Theorem~\ref{abgroup}, we can readily compute all of these
saturated sets; the details are left to the reader. The natural
complementary question concerns saturated sets not containing 0. We can
write $G$ as:

\[
G=\left(\bigoplus_{i=1}^m\frac{\ZZ}{(2k_i+1)\ZZ}\right)\oplus\left(
\bigoplus_{i=1}^n H_i\right)
\]
where the summands $H_i$ are copies of $\ZZ$ or $\frac{\ZZ}{2^n\ZZ}$.
Let the coordinates of $x$ corresponding to the $H_i$'s be 
denoted by $x_i$; the odd elements of $G$ are those with some $x_i$ odd. 
Then the solutions to:

\[
\sum_{i=1}^n a_i x_i \equiv 1\pmod{2}
\]
form a saturated set, where the $a_i$'s are each equal to 0 or 1, and not
all 0; the partition is simply this solution set (which has $|G|/2$
elements) and its complement.

\begin{conj}
The sets described above are the only saturated sets containing no 
even numbers.
\end{conj}

\section{Density Results in $\NN$}
Using the techniques developed in the previous section, we can attack the
classical case $(S,\cdot)=(\NN,+)$. As previously mentioned, few families
of avoidable subsets of $\NN$ have been constructed.
The Fibonacci sequence 
(\cite{Erdos}) was 
the first sequence shown to be avoidable, and other recursive 
sequences (\cite{Shan}, \cite{Develin}, \cite{Evans}) have been analyzed in 
the context of avoidable sets; also, Chow and Long \cite{Chow} 
demonstrated a remarkable connection between avoidable sets and continued 
fractions. All of these sequences grow exponentially; following this lead,
we establish some results on the densities of avoidable sets in $\NN$. We 
start with an observation.

\begin{obs}\label{evensum}
If $a<b<c$, $a+b>c$, and $a+b+c$ is even, then no avoidable set can 
contain $\{a,b,c\}$.
\end{obs}

This is due to the inability to partition $\{\frac{a+b-c}{2},
\frac{a-b+c}{2}, \frac{-a+b+c}{2}\}$ in this case. We now proceed with a
few apropos definitions of density.

For any set $U$, let $U(n)$ be the number of elements of $U$ less than $n$. We 
define the \textbf{arithmetic density} of a set $U$, $d(U)$, by
$d(U)={\displaystyle\limsup_{n\rightarrow\infty}}\frac{U(n)}{n}$. We 
also define 
$d(U>N)={\displaystyle\limsup_{n\rightarrow\infty}}\frac{U(n)-U(N)}{n-N}$.
We note that $d(U)=d(U>N)$ (the proof is routine and will be omitted). 
For a more pertinent measure of density, we define the \textbf{even 
logarithmic density}, $ELD(U)$, and \textbf{logarithmic density}, 
$LD(U)$, as follows, where, $U_2(n)$ is the number of even elements of $U$ 
less than $n$:

\begin{eqnarray*}
ELD(U) &=& \limsup_{n\rightarrow\infty} e^{-\frac{\log 
\frac{n}{2}}{U_2(n)}} \\
LD(U) &=& \limsup_{n\rightarrow\infty} e^{-\frac{\log n}{U(n)}} \\
\end{eqnarray*}

Roughly speaking, $LD(U)$ (resp. $ELD(U)$) is the number $x$ for which 
the 
size of the $n^{th}$ element (resp. $n^{th}$ even element) of $U$ grows 
like $\left(\frac{1}{x}\right)^n$. For instance,
$ELD(\{2^n\})=\frac{1}{2}$, 
$ELD(\{F_n\})=(\frac{1+\sqrt{5}}{2})^{-3}$, where $F_n$ is the 
Fibonacci sequence, and $ELD\left(\left\{\bin{n}{2}\right\}\right)=1$.

Given our past experience, as well as the existence 
of several saturated sets containing only finitely many even 
numbers (such as $\{2,4\}\cup \{x\,|\,x\equiv 3\pmod{8}\}$), it is natural
to consider
avoidable sets with infinitely many even numbers, or to investigate the 
even numbers in a given avoidable set. In this vein, we have
the following two results.
\begin{thm}\label{0density}
Every avoidable set $U$ containing infinitely many even numbers has 
arithmetic density 0.
\end{thm}

\begin{proof}
We want to prove $d(U)<\varepsilon$ for all $\varepsilon>0$. Because 
$d(U)=d(U>N)$, it suffices to find an $N$ for which 
$d(U>N)<\varepsilon$. Let $N$ be an even integer in $U$ 
greater than $\frac{2}{\varepsilon}$, and consider the set $S_k=\{kN+1, 
\ldots, (k+1)N\}$ for $k\ge 1$. No two elements of $U$ with the same 
parity can be in $S_k$, due to Observation~\ref{evensum}. Therefore, 
at most 2 elements of $S_k$ can be in $U$. This clearly implies 
$d(U>N)\le \frac{2}{N} \le \frac{2}{2/\varepsilon}=\varepsilon$, as 
desired.
\end{proof}

\begin{thm}
For any avoidable set $U$, $ELD(U)\le \frac{\sqrt{5}-1}{2}$.
\end{thm}

\begin{proof}
Let $x$ and $y$ be the two smallest even elements of $U$. List the 
even elements of $U$ in increasing order: $\{x,y,a_1, a_2, 
\ldots\}$. Because of Observation~\ref{evensum}, we have $a_1>x+y$, 
$a_2>a_1+y$, and $a_k>a_{k-1}+a_{k-2}$ for $k\ge 3$. 
By induction, this implies $a_k>F_k x+F_{k+1}y$ (using the 
convention $F_1=F_2=1$.) Therefore, there are at most $k+1$ even 
elements of $U$ less than $F_k x + F_{k+1} y$. So we have, if 
$F_{k-1}x+F_k y < n$,
 
\[
\frac{\log \frac{n}{2}}{U_2(n)}\ge \frac{\log (F_{k-1}x+F_k y)-\log 
2}{k+1}
\]
Therefore, we have
\begin{equation}\label{limit}
\liminf_{n\rightarrow\infty}\frac{\log \frac{n}{2}}{U_2(n)}\ge 
\liminf_{k\rightarrow\infty}\frac{\log (F_{k-1}x+F_k y)-\log 2}{k+1}
\end{equation}

The right-hand limit can be evaluated: As $k\rightarrow\infty$, it is 
a well-known fact that
$F_k\rightarrow r\left(\frac{1+\sqrt{5}}{2}\right)^{k}$; this implies that
that the right-hand limit of (\ref{limit}) is equal to 
$\log\frac{1+\sqrt{5}}{2}$ and consequently that 

\[
\liminf_{n\rightarrow\infty}\frac{\log \frac{n}{2}}{U_2(n)}\ge 
\log\frac{1+\sqrt{5}}{2}.
\]

Plugging this into the definition yields $ELD(U)=\frac{\sqrt{5}-1}{2}$ 
as desired.
\end{proof}
This naturally suggests the following conjecture.
\begin{conj}
If $U$ is an avoidable set with infinitely many even numbers, 
$LD(U)\le\frac{\sqrt{5}-1}{2}$.
\end{conj}

The case of $\NN$ illustrates how categorizing avoidable sets becomes 
significantly more difficult when one does not have inverses in 
$S$. In the next section, we discuss another possible difficulty: 
namely, when the group operation $\cdot$ is not commutative.

\section{Families of Non-Abelian Groups}
Note: In this section, all groups will be written multiplicatively, 
and the identity element will be denoted by $e$.

In general, the non-abelian case is more difficult due to the lack of a
convenient description of the group. In the case where the group has a
presentation with few generators and relations, however, we can compute
all of the saturated sets. In this section, we present an example of such
a computation in the case of the dihedral group $D_n$, and the results of
similar computations for other families of groups; the details of these
latter computations are left to the reader, being quite similar to those
of $D_n$. We note that by the definition presented in the introduction, if
$x$ and $y$ are in the same set of a partition avoiding $U\subset G$, we
must have both $xy$ and $yx$ not in $U$.

\subsection{Dihedral Groups}\label{dihed}
Note: In this section, all equalities and variables are in $\Zn$ (or 
$\ZZ/(2n\ZZ)$ when we are investigating $D_{2n}$)
unless noted otherwise, and we use $(a,b)$ to represent the greatest 
common divisor of $a$ and $b$.

$\mathbf{Definition:}$ The $n^{th}$
dihedral group, $D_n$, is defined for $n\ge 3$ to be the group of order $2n$ with
generators $r$ and $f$ satisfying relations $r^n=f^2=1$ and
$fr=r^{-1}f$. It is also the group of symmetries of a regular $n$-gon,
where $r$ represents a rotation by $\frac{2\pi}{n}$ and $f$ represents a
reflection across some axis of symmetry.

We can express elements in $D_n$ uniquely as either $r^a$ or $fr^b$. We
then compute the product law of $D_n$ from the relations between the
generators to be:
\begin{eqnarray*}
r^a \cdot r^b &=& r^{a+n} \\
r^a \cdot fr^b &=& fr^{b-a}\; or\; fr^{b+a} \\
fr^a \cdot fr^b &=& r^{a-b}\; or\; r^{b-a}. \\
\end{eqnarray*}

We will occasionally exploit the definition of $D_n$ as the group of 
symmetries of a regular $n$-gon by labeling an arbitrary element of order 2 $f$, 
instead of starting with a fixed element labeled $f$. Our strategy will be 
to prove a series of partial results, which will coalesce together 
into a characterization of all saturated sets in $D_n$.

\begin{prop}\label{52}
In $D_n$, the saturated sets $U$ containing the identity element are 
$\{e,r^k\}$, where $\frac{n}{(n,k)}$ is even.
\end{prop}

\begin{proof}
First of all, we note that no reflection (element of order 2) can be in an 
avoidable set 
containing $e$. Indeed, label this reflection $f$; then $r, r^{-1}, fr$ is easily 
seen to be a 3-cycle in $G_{\{e,f\}}$. So our avoidable set $U$ must 
consist solely of 
powers of $r$. The connected components of $\GU$ are therefore all either 
contained in $\{0,r,\ldots,r^{n-1}\}$ or disjoint from it. We 
show next that no three-element set $\{e,r^k,r^l\}$ is 
avoidable. Now, 2-coloring 
the components contained in $\{0,r,\ldots,r^{n-1}\}$ is the same as 
2-coloring $G_{\{0,k,l\}}$ in $\Zn$; for this to be 
possible, from Theorem~\ref{cyclic} we know that we must have $2k=l$. 
But 
then $f,fr^k,fr^{2k}$ is a 3-cycle in $\GU$, so $U$ 
is not avoidable. All that remains is to discern when a two-element 
set $\{e,r^k\}$ is avoidable. From Theorem~\ref{cyclic}, we know that 
we can 2-color the components in $\{0,r,\ldots,r^{n-1}\}$. The other 
components are easily seen to be cycles (as $fr^a\cdot fr^a=e$, 
$fr^a$ in $\GU$ has no edges corresponding to $e$) of the form $fr^a, 
fr^{a+k}, fr^{a+2k}, \ldots, fr^{a-k}$, which have length 
$\frac{n}{(n,k)}$. So $U$ is avoidable if and only if this length is 
even, as desired.
\end{proof}
\begin{cor}\label{evencond}
There are no avoidable sets containing $r^k$ if $\frac{n}{(n,k)}$ is odd. 
\end{cor}
Armed with Proposition~\ref{52}, and the fact that
$\{f,fr,\ldots,fr^{n-1}\}$ is saturated (the partition being the set 
and its complement), we conclude the following.
\begin{prop}\label{53}
In $D_{2n+1}$, the saturated sets are exactly $\{e\}$ and $\{f, fr,
\ldots, 
fr^{n-1}\}.$
\end{prop}

We now restrict our attention to $D_{2n}$. As in the case of $D_{2n=1}$, 
$\{f, fr, \ldots, fr^{n-1}\}$ is saturated. Other saturated 
sets with $n$ elements (in each case, the partition consists of the 
set and its complement) are $A\cup B$, $A\cup C$, and $B\cup C$, where
\begin{eqnarray*}
A &=& \{r, r^3, r^5, \ldots, r^{2n-1}\} \\
B &=& \{f, fr^2, fr^4, \ldots, fr^{n-2}\} \\
C &=& \{fr, fr^3, fr^5, \ldots, fr^{n-1}\} \\
\end{eqnarray*}
\begin{prop}\label{54}
In $D_{2n}$, the only saturated sets containing only odd powers of 
$r$ are the three sets described above.
\end{prop}
\begin{proof}
We need only show that no two reflections differing by an odd 
power of $r$ can be in an avoidable set $U$ which contains some
$r^{2k+1}$. 
Label
one reflection $f$, whence the other is $fr^{2l+1}$. Then
these two reflections cannot both be in $U$, for if they were, 
$r^{k-l}, r^{l+k+1}, fr^{l-k}$ would be a 3-cycle in $\GU$.
\end{proof}

So the only uncharacterized avoidable sets are those containing even 
powers of $r$, and not containing $e$. We have the following proposition.

\begin{prop}\label{545}
The only saturated sets containing even powers of $r$ and not $e$ are 
$\{r^{2k}, r^{2l}\}$ where
$\frac{2n}{(2n,2k)}$ and $\frac{2n}{(2n,2l)}$ are even,
and 2 appears with the same 
multiplicity in 2k and 2l.
\end{prop}

\begin{proof}
We assume $U$ contains an even power of $r$, say
$r^{2k}$. $U$ cannot contain any element of the form $r^{2l+1}$,
for then $\GU$ would contain the $(2k+2l+1)$-cycle
\[f, fr^{2k},
\ldots, fr^{(2k)(2l)}, fr^{(2k)(2l+1)}, fr^{(2k-1)(2l+1)},
\ldots, fr^{2l+1}.\]
Furthermore, $U$ cannot contain a
reflection $f$, for then $G_U$ contains the 3-cycle $fr^{n+k}, fr^{n-k},
r^{n-k}$.
So if $U$ is to have any other elements, they must be even
powers of $r$. There cannot be two other even powers of $r$, because
then the set would not even be avoidable in the subgroup
$\ZZ/(2n\ZZ)$ by Theorem~\ref{cyclic}. Consequently, we only need to
discern when $U=\{r^{2k}, r^{2l}\}$ is avoidable. We know from
Corollary~\ref{evencond} that $\frac{2n}{(2n,2k)}$ and
$\frac{2n}{(2n,2l)}$ must be even, and we know that the components of
$\GU$ contained in $\{e, r, r^2, \ldots, r^{2n-1}\}$ are 2-colorable
from Theorem~\ref{cyclic}. So we need only color the component 
containing an arbitrary reflection $f$.

Looking at the associated graph, we can easily see using the product
laws for $D_n$ that the elements in the connected component of $f$ 
are those of the form $fr^{a(2k)+b(2l)}$ with this point being an odd 
distance from $f$ if and only if $a+b$ is odd. So $\GU$ will contain an
odd cycle if and only if the equation
\begin{equation}\label{critical}
a(2k)+b(2l)=0\pmod{2n}
\end{equation}
has any solutions with $a+b$ odd. Note that as $2k$ and $2l$ are both 
nonzero, and $\frac{2n}{(2n,2k)}$ and $\frac{2n}{(2n,2l)}$ are even,
this equation has no trivial solutions (i.e. with one of 
$a$, $b$ equal to 0.)
\begin{claim}
Equation~(\ref{critical}) has solutions with $a+b$ odd if and only if $2k$ and 
$2l$ have different multiplicities of 2.
\end{claim}
\textit{Proof.}
Suppose $2k$ and $2l$ have different multiplicities of 2. Then 
exactly 
one of $\frac{lcm(2k,2l)}{2k}$ and -$\frac{lcm(2k,2l)}{2l}$ is odd. Letting 
these equal $a$ and $b$, respectively, yields the desired solution.

Conversely, suppose $2k$ and $2l$ have equal multiplicities of 2, say $s$. Then 
2 appears with multiplicity $t>s$ in $2n$, because 
$\frac{2n}{(2n,2k)}$ is even. But now consider $a(2k)+b(2l)$ for 
$a+b$ odd. Without loss of generality, $a$ is odd and $b$ is even. 
But then $2^{s+1}$ divides $b(2l)$ but not $a(2k)$; consequently, it 
does not divide their sum. However, $2^{s+1}$ divides $2^t$, which 
divides $2n$, so $a(2k)+b(2l)$ cannot be equal to 0 (mod $2n$) for 
$a+b$ odd. This proves the claim.

Therefore, since the existence of solutions to 
equation~(\ref{critical}) with $a+b$ odd is equivalent to the 
unavoidability of $\{r^{2k},r^{2l}\}$, Proposition~\ref{545} is proven.
\end{proof}

We now have enough information to characterize all avoidable sets in
$D_n$.
\begin{thm}\label{dnsat}
The saturated sets in $D_n$ are:

(a) if $n$ is odd, $\{e\}$ and $\{f, fr, \ldots, fr^{n-1}\}$

(b) if $n$ is even, $\{e, r^k\}$ where $\frac{n}{(n,k)}$ is even; 
$\{r^{2k}, r^{2l}\}$ where $\frac{n}{(n,2k)}$ and $\frac{n}{(n,2l)}$ 
are even and 2 appears with equal multiplicity in $2k$ and $2l$; and 
the three $n$-element sets $A\cup B$, $A\cup C$ and $B\cup C$, where 
$A$ is the set $\{r, r^3, \ldots, r^{n-1}\}$, $B$ is the set $\{f, fr^2,
\ldots, fr^{n-2}\}$, and $C$ is the set $\{fr, fr^3, \ldots, fr^{n-1}\}$.
\end{thm}

\subsection{Semi-Dihedral Groups}
Closely related to dihedral groups are semi-dihedral groups. 
We define for $m>3$ $SD_m$, the \textbf{$m^{th}$ semi-dihedral group}, 
to be the group with generators $x$, $y$ satisfying relations 
$x^{2^{m-1}}=y^2=e$ and $yx=x^k y$ where $k=-1+2^{m-2}$. Noting that 
$k^2=1$, and writing all elements as $x^a$ or $yx^a$, we can easily 
derive the multiplication law for semi-dihedral groups to be:
\begin{eqnarray*}
x^a \cdot x^b &=& x^{a+b} \\
x^a \cdot yx^b &=& yx^{a+b}\; or\; yx^{a+kb} \\
yx^a \cdot yx^b &=& x^{ka+b}\; or\; x^{kb+a}. \\
\end{eqnarray*}

By methods similar to those employed in Section~\ref{dihed}, we arrive
at the following theorem for semi-dihedral groups. 

\begin{thm}
The saturated sets in $SD_m$ are $\{e,x^{2^{m-2}}\}$;
$\{x^{2n}, x^r\}$, where $2n, r\ne 0, 2^{m-2}$ and $2n$ and $r$ have 
the same divisibility by 2; and $A\cup B$, $A\cup C$, $B\cup C$
where $A=\{x, x^3, \ldots, x^{2^{m-1}-1}\}$, $B=\{y, yx^2, \ldots, 
yx^{2^{m-1}-2}\}$, and $C=\{yx, yx^3, \ldots, yx^{2^{m-1}-1}\}$.
\end{thm}
\subsection{Generalized Quaternion Groups}
In the same category as the dihedral and semi-dihedral groups are the 
generalized quaternion groups. We define $Q_m$, the \textbf{$m^{th}$ 
generalized quaternion group}, to be the group defined with generators 
$a,b$ and relations $a^{4m}=b^2=1$, $a^{2m}=b^2$, and $ba=a^{-1}b$. As 
before, we can represent each element as $a^x$ or $ba^x$, whence the 
multiplication law for $Q_m$ is as follows:

\begin{eqnarray*}
a^x \cdot a^y &=& a^{x+y} \\
ba^x \cdot a^y &=& ba^{x+y}\; or\; ba^{x-y} \\
ba^x \cdot ba^y &=& a^{2m+x-y}\; or\; a^{2m+y-x}. \\
\end{eqnarray*}

Again using the methods exhibited in Section~\ref{dihed}, we can compute
all saturated subsets of $Q_m$.
\begin{thm}\label{614}
The saturated sets in $Q_m$ are:
\begin{itemize}
\item $\{a^{2m},a^{2r}\}$ where $2r$ and $2m$ have different 
multiplicities of 2, or $2r\equiv 0\pmod{4m}$
\item $\{a^{2m},ba^r\}$
\item $\{e,a^{2r}\}$ where $\frac{2r}{(2m,2r)}$ is even
\item $\{a^{2r},a^{2s}\}$ where $2r, 2s\ne 0, 2m \pmod{4m}$; $2r$ and $2s$ 
have the same multiplicity of 2; and this multiplicity is different from that 
of $2m$
\item $A\cup B$, $A\cup C$, and $B\cup C$, where:
\begin{eqnarray*}
A &=& \{a, a^3, \ldots, a^{4m-1}\} \\
B &=& \{b, ba^2, \ldots, ba^{4m-2}\} \\
C &=& \{ba, ba^3, \ldots, ba^{4m-1}\}. \\
\end{eqnarray*}
\end{itemize}
\end{thm}
\subsection{Non-abelian Groups of Order $pq$}
Another family for which we can readily compute the saturated sets is
the class of non-abelian groups of order $pq$. The only
nonabelian groups of order $pq$, where $p$ and $q$ are prime and $p>q$,
occur when $p\equiv 1\pmod{q}$. In this case,
there exists a unique group given by the presentation $G$, generated by 
$a,b$ with relations $a^q=b^p=e$ and $ba=ab^s$, where $s\ne 
1$ and $s^q=1$; we
assume $p,q$ are odd because if $q=2$, $G$ is simply the dihedral group $D_p$. 
We represent all elements of the 
group as $a^x b^y$ with $x\in \frac{\ZZ}{q\ZZ}$ and $y\in \frac{\ZZ}{p\ZZ}$. 
Furthermore, we denote by $\frac{k}{2}$ the element $x$ such that $2x=k$ in 
$\frac{\ZZ}{q\ZZ}$. We present the product law:
\[
a^k b^l \cdot a^m b^n = a^{k+m} b^{ls^m+n}\; or\; a^{k+m}b^{l+ns^k}.
\]

In this case, the algebra is somewhat more subtle than that of $D_n$. To
give a flavor of the methods used, we present the proof of the following
main lemma. 

\begin{lem}\label{73}
For $k\ne 0$, the set $\{a^k b^l, a^m b^n\}$ is not 
avoidable 
unless $a^m b^n$ = $(a^k b^l)^r.$ 
\end{lem}
\begin{proof}
Consider the set $\{a^{\khalf} 
b^x, a^{\khalf} b^y, a^{m-\khalf} b^z\}$. This set is a 3-cycle in the
associated graph of this set if and only if $x\neq y$ and the elements $x,
y,$ and $z$ of the finite field
$\mathbb{F}_p$ satisfy either:
\begin{eqnarray*}
x+s^{\khalf}y &=& l \\
s^{m-\khalf}y + z &=& n \\
x+s^{\khalf}z &=& n \\
\end{eqnarray*}
or
\begin{eqnarray*}
s^{\khalf}x+y &=& l \\
s^{m-\khalf}y + z &=& n \\
x+s^{\khalf}z &=& n. \\
\end{eqnarray*}
Now, the system of equations $A\mathbf{v}=\mathbf{c}$ 
has a solution if $\det A\ne 0$. For the first set of equations, $\det 
A=s^m+s^{\khalf}$; for the second set, $\det A=s^{m+\khalf}+1$. But these 
cannot both be 0; if the second equation is, we obtain $s^m=-s^{-\khalf}$, 
whereupon the determinant associated to the first set of equations is equal to 
$s^{\khalf}-s^{-\khalf}$. This quantity is then nonzero as
$s^{\khalf}=s^{-\khalf}$ implies $s^k=1$, which is impossible as $k\ne
0$. Therefore, one set of equations must 
have a solution; we need to show that this solution has $x\ne y$. If $x=y$ and $(x,y,z)$ satisfies either 
set of equations, we obtain, substituting into either set of equations:
\begin{eqnarray*}
x(1+s^{\khalf}) &=& l \\
s^{m-\khalf}x+z &=& n \\
x+s^{\khalf} z &=& n \\
\end{eqnarray*}
Multiplying the second equation by $s^{\khalf}$ and subtracting the third yields
\[
(s^m-1)x=n(s^{\khalf}-1).
\]
Now, we multiply on both sides by $(1+s^{\khalf})$ to get
\begin{equation}\label{klmn}
(s^m-1)l=n(s^k-1).
\end{equation}
However, let $\frac{m}{k}$ be the number (mod $q$) which when multiplied by $k$ 
(mod $q$) yields $m$; we then have:
\begin{eqnarray*}
(a^k b^l)^{\frac{m}{k}} &=& a^m b^{l+ls^k+\cdots+ls^{m-k}} \\
&=& a^m b^{l(\frac{s^m-1}{s^k-1})} \\
&=& a^m b^n \\
\end{eqnarray*}
by substitution from equation~(\ref{klmn}), so $a^m b^n=(a^k b^l)^r$ for 
$r=\frac{m}{k}$. Consequently, unless this condition holds, we can find a
bona fide 
3-cycle in $G_{\{a^k b^l, a^m b^n\}}$, so Lemma~\ref{73} is proven.
\end{proof}

This lemma is then the main step in the proof of the following
surprisingly simple theorem categorizing all saturated sets in a
non-abelian group of order $pq$.
\begin{thm}
The saturated sets in a non-abelian group of order $pq$ are all sets of the form 
$\{x, x^n\}$.
\end{thm}

\section{Conclusion}

In general, the case of non-abelian groups is harder than the case of
abelian groups, because the group may not have a simple,
easy-to-manipulate, presentation. There are some recurring themes,
however: all of the saturated sets in the cases of both abelian and
non-abelian groups are either small relative to the order of the group (in
the families we have discussed, constant with respect to the order of the
group) or the nontrivial coset of a subgroup of order 2. It is in general
true that the nontrivial coset of a subgroup of order 2 is avoidable (and
saturated unless all elements of the subgroup have order 2); an
interesting open question is whether or not all other avoidable sets are
small relative to the size of the group. 

Other desirable results would include linking the avoidable sets in a
group to the avoidable sets in a subgroup or quotient group. In the first case, the problem 
with passing to a subgroup is that saturated sets entirely contained in the 
subgroup will of course still be avoidable, but may not be saturated; for 
example, in the group $D_{2n+1}$, the subset $\{e\}$ is saturated, but it is 
certainly not saturated in the subgroup $\frac{\ZZ}{(2n+1)\ZZ}$. In the second 
case, we would like to lift avoidable sets in the quotient group to avoidable 
sets in the larger group simply by lifting the partition; however, two distinct 
lifts of a single element in the quotient group may multiply to some element
in the lift of the avoidable set.

If these results and a characterization of avoidable sets in the symmetric
group $S_n$ were obtained, the problem for finitely generated groups would be
solved. However, the latter also proves daunting. In $S_n$, the set of all odd
permutations is always avoidable (and saturated for $n\ge 3$); in $S_4$, all
the other saturated sets have size 1. It is not true in general that all the
other saturated sets have size 1, though; for instance, in $S_6$, the set
$\{(12)(34), (56)\}$ is avoidable. Nevertheless, saturated subsets of $S_n$
other than the set of odd permutations tend to be quite small.

\section*{Acknowledgements}
This work was done under the supervision of Joseph Gallian at the University of
Minnesota, Duluth, with financial support from the National Science Foundation
(Grant DMS-9531373-001) and the National Security Agency (Grant
MDA-904-98-1-0523). The author wishes to thank Joseph Gallian and Daniel Biss
for their helpful comments on preliminary versions on this manuscript, and
Daniel Isaksen for his suggestions during the research presented in this
paper.

\end{document}